\documentclass[12pt,a4paper,draft]{article}
\usepackage[english]{babel}
\usepackage{latexsym,amsfonts,amssymb,amsmath,longtable,amsthm,mathrsfs,cite}
\renewcommand{\le}{\leqslant}

\renewcommand{\qed}{\hfill{$\Box$}}

\newtheorem*{Prop*}{{\bfseries Proposition}}\newtheorem*{Cor*}{{\bfseries Corollary}}

\newtheorem{Lemma}{{\bfseries Lemma}}
\newtheorem{Cor}{{\bfseries Corollary}}
\newtheorem{Theo}{{\bfseries Theorem}}
\theoremstyle{definition}

\newtheorem{Rem}{{\bfseries Remark}}

\DeclareMathOperator{\GL}{GL} \DeclareMathOperator{\Aut}{Aut}
\DeclareMathOperator{\Sym}{Sym}
\DeclareMathOperator{\SL}{SL}\DeclareMathOperator{\Hall}{Hall}

\title{\vspace{-1cm} \hfill{\normalsize 20D20, 20D30}{
\fontfamily{cmr} \fontseries{bx} \selectfont \\ \vspace{1cm} Frattini Argument for Hall subgroups}
\thanks{The work is supported by  RFBR, projects  11-01-00456 and
12-01-33102.}}
\date{}
\author{\bf  D.O.Revin, E.P.Vdovin}

\begin{document}
\sloppy


\maketitle
\pagenumbering{arabic}
\begin{abstract}

In the paper, it is proved that if a finite group $G$ possesses a $\pi$-Hall subgroup for a set $\pi$ of primes, then every normal subgroup
$A$ of $G$ possesses a $\pi$-Hall subgroup $H$ such that ${G=AN_G(H)}$.

\end{abstract}


\section*{Introduction}

One of the most frequently used statement in the finite group theory is the following evident corollary to the Sylow theorem on the conjugacy of Sylow  $p$-subgroups \cite{Frattini}, \cite[(6.2)]{Aschbacher}.

\begin{Prop*} {\rm \bf(Frattini Argument)}
Let $A$ be a normal subgroup of a finite group $G$,  and let $S$ be a Sylow $p$-subgroup of $A$ for a prime $p$. Then ${G=AN_G(S)}.$
\end{Prop*}

Notice that if $A$ is normal in $G$, then the action of $G$ by conjugation on the set of subgroups of $A$ induces the action of $G$ on
the set $$\{H^A\mid H\leq A \}$$ of conjugacy classes of subgroups of $A$ (here and thereafter, for a subgroup $H$ of $G$, we set
$H^G=\{H^g\mid g\in G\}$). Moreover, the stabilizer of the class $H^A$ consisting of subgroups conjugate with $H$ in $A$ coincides with $AN_G(H)$.
Thus, the equality $G=AN_G(H)$ is equivalent
for  $H^A$ to be invariant under the action of $G$, or in other words, is equivalent to the equality~${H^A=H^G}$  (cf. \cite[(6.3)]{Aschbacher} for example).

Analogs of the Frattini Argument are valid not only for Sylow subgroups, but also for some other classes of subgroups. For example, it is proved in \cite{ZMR} that every normal
subgroup   $A$ of a finite group  $G$ possesses a maximal solvable subgroup $S$ such that  $G=AN_G(S)$, and, as a consequence, in every finite group, there is a subgroup that is a solvable injector
and a solvable projector simultaneously.

The notion of a $\pi$-Hall subgroup is a natural generalization of the notion of Sylow subgroup.
Let $\pi$ be a set of primes. A subgroup   $H$ of a finite group $G$ is called a {\it $\pi$-Hall subgroup},
if all prime dividers of the order of  $H$ lies in $\pi$, while the index of $H$ is not divisible by the elements of~$\pi$. We denote by $\Hall_\pi(G)$ the set of $\pi$-Hall subgroups
of a finite group $G$.

The Hall theorem states that, in a solvable group,   $\pi$-Hall subgroups exist and are conjugate for every set $\pi$ of primes. As a corollary, in every finite solvable group, for every set  $\pi$ of primes, a complete analogue of the Frattini argument for $\pi$-Hall subgroups holds. 
However, given a set   $\pi$ of primes, in a nonsolvable finite group $\pi$-Hall subgroups may fail to exist, and, even if they exist, they can be nonconjugate.

According to  \cite{Hall}, we say that  $G$ {\it satisfies  $E_\pi$} (or, briefly  $G\in E_\pi$) if $\Hall_\pi(G)\ne\varnothing$. If $G\in E_\pi$ and all  $\pi$-Hall subgroups are conjugate then we say that   $G$ {\it satisfies
$C_\pi$} ($G\in C_\pi$). A group satisfying $E_\pi$ or $C_\pi$ is called also an $E_\pi$- or
a $C_\pi$-{\it group}, respectively. 

It is easy to show that if  $A$ is a normal subgroup of a finite group $G$ and $A\in C_\pi$ for some  $\pi$ then  $G=AN_G(H)$ for every $\pi$-Hall
subgroup $H$ of $A$.
This notice is a key idea in the proof of the famous Chunikhin Theorem \cite[(3.12)]{Suz}, claiming that an extension of a $C_\pi$-group by a  $C_\pi$-group is a  $C_\pi$-group as well.
In general, distinct analogues of the Frattini Argument for Hall subgroup are of extraordinary importance in the study of Hall subgroups and generalizations of the Sylow Theorem
(see \cite{Hall,Chun1,Chun2,Chun3,Wie,Wie2,Wie3,Wie4,Hart,Shem3,Shem4e,GrossExistence,RevVdoArxive,ExCrit,MR,R,VR,VR1,VR2,VR3,VR4,Baer}).

The following statement is the main result of the paper.

\begin{Theo}\label{MainTheorem} {\em (Frattini Argument for Hall subgroups)}
Let $G\in E_\pi$ for some set $\pi$ of primes and $A\unlhd G$. Then there exists a $\pi$-Hall subgroup  $H$ of $A$ such that $G=AN_G(H)$. Moreover  $N_G(H)\in E_\pi$ and every  $\pi$-Hall subgroup of   $N_G(H)$ is a $\pi$-Hall subgroup of~$G$.
\end{Theo}

\begin{Rem}\label{R1} It is well known that, if $A\unlhd G$ and $H$ is a $\pi$-Hall subgroup of~$G$, then  $H\cap A$ is a $\pi$-Hall subgroup
of $A$ (see Lemma \ref{base}). Thus, the condition  $G\in E_\pi$ in Theorem~\ref{MainTheorem} implies $A\in E_\pi$. However we cannot replace the condition $G\in E_\pi$ with a
weaker condition $A\in E_\pi$. Indeed, let  $\pi=\{2,3\}$ and $A=\GL_3(2)=\SL_3(2)$. Then $A$ possesses exactly two classes of conjugate $\pi$-Hall subgroups with the representatives
$$
H_1=\left(
\begin{array}{c@{}c}
\fbox{
$\begin{array}{c}
\!\!\!
\GL_2(2)
\!\!\! \\
\end{array}$
}
& *\\
0 &\fbox{1}
\end{array}
\right){\text{\,\,\, и\,\,\, }}
H_2=\left(
\begin{array}{c@{}c}
\fbox{1}& *\\
 0&\fbox{$\begin{array}{c}
\!
\GL_2(2)
\!
\\
\end{array}$}
\end{array}
\right).
$$
 The class $H_1^A$ consists of the line stabilizers in the natural module for
$G$, while  $H_2^A$ consists of the plane stabilizers. The map  $\iota : x\in A\mapsto (x^t)^{-1}$ is an automorphism
of order  2 of $A$ (here and below $x^t$ is the transpose of  $x$). This automorphism interchanges the classes $H_1^A$ and $H_2^A$. Consider  the natural split extension
$G=A:\langle\iota\rangle$. Then $N_G(H)\leq A$ for every $\pi$-Hall subgroup  $H$ of  $A$ and, in particular,  $AN_G(H)=A<G$.  Indeed, since   $|G:A|=2$, the assertion 
$N_G(H)\nleq A$ would imply $G=AN_G(H)$ and $H^G=H^A$,  in particular $H^\iota\in H^A$. But in this case the automorphism  $\iota$ would leave the class $H^A$   invariant, a
contradiction. Now  $G\notin E_\pi$, since otherwise, in view of the identity $|G:A|=2$, the equality $G=AK$ would hold for every $\pi$-Hall subgroup  $K$ of $G$ and, in
particular, $K\nleq A$. Since $K\cap A$ is a $\pi$-Hall subgroup of  $A$ and $K\leq N_G(K\cap A)$ it follows that $G=AN_G(K\cap A)$, a contradiction.
\end{Rem}

\begin{Rem} If one compare the conclusion of Theorem~\ref{MainTheorem} with the original Frattini Argument, it will be natural natural to ask, whether the identity $G=AN_G(H)$ holds for
every $\pi$-Hall subgroup $H$ of $A$, or at least for every subgroup of the form $H=K\cap A$, where $K$ is a $\pi$-Hall subgroup of $G$, if the condition of Theorem~\ref{MainTheorem} is satisfied? This question has a negative answer. Indeed, let $\pi=\{2,3\}$ and $S=\GL_3(2)$. Suppose that subgroups $H_1$ and $H_2$ are defined in the same way as in Remark \ref{R1}. Consider the direct product
$$
A=\underbrace{S\times\dots\times S}\limits_{5\text{ times}}.
$$
Now $A$ admits
an automorphism $\tau:(x_1,x_2,\dots,x_5)\mapsto (x_5,x_1,\dots,x_4)$ and we denote by  $G$ the natural semidirect product of $A$ and~$\langle \tau\rangle$. Since $|G:A|=5$ is a prime, $5\notin\pi$, and  $A\unlhd G$,  the sets  $\Hall_\pi(G)$ and $\Hall_\pi(A)$ coincide. Consider $\pi$-Hall subgroups $$K_1=\underbrace{H_1\times\dots\times H_1}\limits_{4\text{ times}}\times H_2,$$ $$ K_2=H_2  \times\underbrace{H_1\times\dots\times H_1}\limits_{4\text{ times}}
$$
of $A$ (hence also  of  $G$). Clearly,  $K_1=K_2^\tau$, but $K_1$  $K_2$ are not conjugate in  $A$. Thus, the conjugacy classes of  $K_1$ and $K_2$ in $A$ are distinct,  while they are fused in $G$. Whence $G\ne AN_G(K_i)$, $i=1,2$.
\end{Rem}

As a corollary to Theorem~\ref{MainTheorem} we obtain the main result of \cite{VR1}, thus we give a simpler proof to the following statement.
\begin{Cor} \label{HN}
If $G\in C_\pi$, $A\unlhd G$ and $H$ is a $\pi$-Hall subgroup of $G$, then $HA\in C_\pi$.
\end{Cor}

\begin{Cor} \label{E_piCrit}
Let  $A$ be a normal subgroup of $G$. Then  $G\in E_\pi$ if and only if $A\in E_\pi$, $G/A\in E_\pi$, and there exists $H\in \Hall_\pi(A)$ such that $H^A=H^G$.
\end{Cor}

\begin{Cor} \label{ExInvHallSubg}
Let   $G\in E_\pi$, $A\leq\Aut(G)$ and $(|G|,|A|)=1$. Then there exists an $A$-invariant $\pi$-Hall subgroup $H$ of $G$.
\end{Cor}
\section{Preliminary results}

We always denote by $\pi$ a set of primes and by $\pi'$ the complement to $\pi$ in the set of all primes. An integer $n$ is called a $\pi$-number if every prime divisor of $n$ belongs to $\pi$. A group $G$ is called a $\pi$-group if $|G|$ is a $\pi$-number.

\begin{Lemma} \label{base}
Let $A$ be a normal subgroup of  $G$.  If $H$ is a $\pi$-Hall subgroup of $G$, then
$H \cap A$ is a $\pi$-Hall subgroup of $A$, and $HA/A$ is a $\pi$-Hall subgroup of~$G/A$.
\end{Lemma}

\begin{proof} See \cite[Ch. IV, (5.11)]{Suz}).
\end{proof}

Recall that a finite group $G$ is said to be {\it $\pi$-separable}, if there is  a normal series of $G$ with all  factors being either  $\pi$- or $\pi'$-groups.

\begin{Lemma} \label{base1}
Every $\pi$-separable group satisfies~$C_\pi$.
\end{Lemma}
\begin{proof} See \cite[Ch. V, Theorem~3.7]{Suz}.
\end{proof}

\begin{Lemma}\label{HallExist}
Let $A$ be a normal subgroup of $G$ such that $G/A$ is a $\pi$-group, $U$ be a $\pi$-Hall subgroup of $A$. Then a
$\pi$-Hall subgroup $H$ of $G$ with $H\cap A=U$ exists if and only if $U^G=U^A$.
\end{Lemma}
\begin{proof} See \cite[Lemma~2.1(e)]{RevVdoArxive}.
\end{proof}

\begin{Lemma} \label{HallSubgrIn HomImage}
Let $A$ be a normal subgroup of an $E_\pi$-group $G$. Then for every $K/A\in\Hall_\pi(G/A)$ there exists $H\in\Hall_\pi(G)$ such  that $K=HA$.
\end{Lemma}

\begin{proof} See \cite[Corollary 9]{ExCrit}.
\end{proof}

If $S$ is a subnormal subgroup of $G$, then by $\Hall_\pi^G(S)$ we denote the set of subgroups of type   $H\cap S$, where $H\in\Hall_\pi(G)$. Clearly $\Hall^G_\pi(S)$ is a union of several classes of conjugate $\pi$-Hall subgroups of $S$, and let $k^G_\pi(S)$ be the number of these classes.

Recall that a subgroup generated by all minimal normal subgroups is called the {\it socle}. A group is called {\it almost simple}, if its socle is a finite simple nonabelian group.

\begin{Lemma} \label{Simple}
Let $\pi$ be a set of primes and $G$ be an almost simple $E_\pi$-group with $($nonabelian simple$)$ socle~$S$. Then the following hold.

\begin{itemize}
 \item[$(1)$] If $2\not\in\pi$, then  $k_\pi^G(S)=1.$

\item[$(2)$] If $3\not\in\pi$, then $k_\pi^G(S)\in\{1,2\}$.

\item[$(3)$] If $2,3\in\pi$, then $k_\pi^G(S)\in\{1,2,3,4,9\}$.
\end{itemize}
In particular, $k_\pi^G(S)$ is a $\pi$-number.
\end{Lemma}

\begin{proof} See \cite[Theorem~1.1]{RevVdoArxive}.
\end{proof}

\begin{Lemma} \label{MainLemma}
Let $\pi$ be a set of primes, $G$ be an almost simple  $E_\pi$-group with socle~$S$ and $T\leq G$. Consider the action of $T$ on $$\Omega=\{(H\cap S)^S\mid H\in\Hall_\pi(G)\}$$ by conjugation. Then the following hold:
\begin{itemize}
\item[$(1)$] if $k^G_\pi(S)\ne 9$, then the length of every orbit is a $\pi$-number;
\item[$(2)$] if $k^G_\pi(S)= 9$, then there exists an orbit such that its length is a $\pi$-number.
\end{itemize}
\end{Lemma}
\begin{proof}
By definition, $k^G_\pi(S)=|\Omega|$.  Consider an orbit $\Delta$ of $T$ on $\Omega$. Since   $|\Delta|\le |\Omega|$, Lemma
\ref{Simple} implies that either  $|\Delta|$ is a
$\pi$-number, or   $2,3\in\pi$, $k_\pi^G(S)= 9$, and $|\Delta|\in \{5,7\}$. In the last case, the action of  $T$ on $\Omega$ is intransitive,  every orbit under this action with the exception of  $\Delta$ has the length at most $|\Omega|-|\Delta|\le 4$, and so is a $\pi$-number.
\end{proof}

\begin{Lemma} \label{MainLemma1}
Let $G$ be an almost simple $E_\pi$-group, $S$ be the socle of, and $H\in\Hall_\pi(G)$. Then, for some $H\in\Hall_\pi(G)$, the equality $(H\cap S)^G=(H\cap S)^S$ holds.
\end{Lemma}

\begin{proof} Consider $$\Omega=\{(H\cap S)^S\mid H\in \Hall_\pi(G)\}.$$
Let $T/S$ be a $\pi'$-Hall subgroup of a (solvable) group $G/S$. Consider the action of $T$ on $\Omega$ induced by the action of $T$ on the set of subgroups of $S$ via conjugation. By Lemma \ref{MainLemma} it follows that $T$ possesses an orbit~$\Delta$ of length being a $\pi$-number. On the other hand, since $S$  is included in the kernel of the action of $T$ on $\Omega$ and since  $T/S$ is a  $\pi'$-group we obtain that $|\Delta|$ is a  $\pi'$-number. Hence  $|\Delta|=1$.

Thus,  $(H\cap S)^T=(H\cap S)^S$ for some
$H\in\Hall_\pi(G)$.  Since $G=HT$ we also have $$(H\cap S)^G=(H\cap S)^{HT}=(H\cap S)^T=(H\cap S)^S,$$
and the lemma follows.
\end{proof}

Let $A,B,H$ be subgroups of  $G$ such that $B\unlhd A$. By $N_H(A/B)$ we denote the intersection  $N_H(A)\cap N_H(B)$. Then each element $x\in
N_H(A/B)$ induces an automorphism on $A/B$ acting by $Ba\mapsto B x^{-1}ax$. Thus there exists a homomorphism
$N_H(A/B)\rightarrow \Aut(A/B)$. Denote   by $\Aut_H(A/B)$ the image of the homomorphism and call it  a {\em  group of
$H$-induced automorphisms on $A/B$}. The kernel of the homomorphism is denoted by~${C_H(A/B)}$. If  $B=1$, then  $\Aut_H(A/B)$ is denoted by~$\Aut_H(A)$.

\begin{Lemma} \label{InvGross}
Let $S$ be a simple subnormal subgroup of an $E_\pi$-group $G$. Then $\Aut_G(S)\in E_\pi$.
\end{Lemma}

\begin{proof} There exists a composition series
$$1=G_0<G_1\dots<G_n=G$$ of $G$ such that $S=G_1$. The claim of the Lemma follows from \cite[Theorem~4]{ExCrit} and the equality $\Aut_G(G_1/G_0)=\Aut_G(S)$.
\end{proof}

\begin{Lemma} \label{HallSubgrInMinNorm}
Let $S$ be a simple subnormal  $E_\pi$-subgroup of $G$ and
$U\in\Hall_\pi(S)$. Choose a right transversal $g_1,\dots,g_n$ for   $N_G(S)$ in $G$ and let $$V=
\langle U^{g_i}\mid i=1,\dots, n\rangle.$$ Then $V\in\Hall_\pi(\langle
S^G\rangle)$.
\end{Lemma}

\begin{proof} One may assume that $S$ is nonabelian. Let $A=\langle S^G\rangle$. Since $S$ is nonabelian simple, we derive  that $A$ is a direct product of  $S^{g_1},\dots, S^{g_n}$ and, as a consequence,  $V$ is a direct product of   $U^{g_1},\dots, U^{g_n}$. In particular,
$V$ is a  $\pi$-group and $|A:V|=|S:U|^n$ is a  $\pi'$-number, i.e.
$V\in\Hall_\pi(A)$.
\end{proof}

\begin{Lemma} \label{InvarClass}
Let $S$ be a simple subnormal subgroup of $G$ and a subgroup $U$ of $S$ is chosen so that the equality
$U^{\Aut_G(S)}=U^S$ holds. Take a right transversal  $g_1,\dots,g_n$ for $N_G(S)$ in $G$. Let $$A=\langle S^{g_i}\mid
i=1,\dots, n\rangle\text{ and }V= \langle U^{g_i}\mid i=1,\dots, n\rangle.$$ Then
$A=\langle S^G \rangle\unlhd G$ and $V^G=V^A$. Moreover, if $U\in\Hall_\pi(S)$, then $V\in\Hall_\pi(A)$.
\end{Lemma}

\begin{proof} In view of the choice of elements $g_1,\dots,g_n$, for every $g\in G$ we have
$S^g\in\{S^{g_i}\mid i=1,\dots, n\}$ and so $A=\langle S^G \rangle\unlhd G$.
Notice also that since  $S$ is simple and subnormal, we have
$[S^{g_i},S^{g_j}]=1$ for $i\ne j$.

Let $g\in G$. There exist a permutation  $\sigma\in\Sym_n$  and $x_1,\dots,x_n\in N_G(S)$ such that
$g_ig=x_ig_{i\sigma}$. Now consider $\gamma_i\in\Aut_G(S)$, where the map
$\gamma_i:S\rightarrow S$ is given by $s\mapsto s^{x_i}$. By condition,
$U^{x_i}=U^{\gamma_i}=U^{s_i}$ for some $s_i\in S$. Set
$a_i=s_{i\sigma^{-1}}^{g_i}$ and $a=a_1\cdot \ldots\cdot a_n$. Clearly $a\in A$. We remain to show that  $V^g=V^a$, and therefore the equality $V^G=V^A$ holds.

By the definition, it follows that
$a_i\in S^{g_i}$ and $U^{g_ia}=U^{g_ia_i}$. We have
\begin{multline*}
V^g= \langle U^{g_ig}\mid i=1,\dots, n\rangle= \langle U^{x_ig_{i\sigma}}\mid
i=1,\dots, n\rangle=\\ \langle U^{s_ig_{i\sigma}}\mid i=1,\dots, n\rangle=
 \langle U^{s_{i\sigma^{-1}}g_{i}}\mid i=1,\dots, n\rangle=\\ \langle
U^{g_is_{i\sigma^{-1}}^{g_i}}\mid i=1,\dots, n\rangle=\langle U^{g_ia_i}\mid
i=1,\dots, n\rangle=\\
 \langle U^{g_ia}\mid i=1,\dots, n\rangle=V^a,
\end{multline*}
and the lemma follows.
\end{proof}

\section{Proof of the main results}

{\it Proof of Theorem}~\ref{MainTheorem}. Let $G\in E_\pi$ and $A\unlhd G$. We show by induction by the order of $G$, that $A$ possesses
a $\pi$-Hall subgroup $H$ such that $G=AN_G(H)$.

If $G$ is simple, we have nothing to prove.

Firstly, consider the case, where   $A$ is a minimal normal subgroup of $G$. Let $S$ be a minimal subnormal subgroup of  $A$. By Lemma \ref{InvGross} we have $\Aut_G(S)\in E_\pi$. Lemma \ref{MainLemma1} implies that $S$ possesses a
$\pi$-Hall subgroup $U$ such that $U^{\Aut_G(S)}=U^S$. Choose a right transversal $g_1,\dots, g_n$ of $G$ by
$N_G(S)$, and let $$H=\langle U^g_i\mid i=1,\dots, n\rangle.$$ Since $A$ is a minimal normal subgroup of $G$ we obtain $A=\langle
S^G\rangle$. By Lemma \ref{HallSubgrInMinNorm}, $H$ is a $\pi$-Hall subgroup of $A$. By Lemma
\ref{InvarClass} it follows that $H^G=H^A$ or equivalently $G=AN_G(H)$.

Thus we may assume that there exists a minimal normal subgroup  $M$ of $G$ such that $M<A$. As we proved above, the identity $G=MN_G(V)$ holds for a  $\pi$-Hall subgroup  $V$ of $M$. Set  $K=N_G(V)$. Notice that by Lemma \ref{HallExist} the identity $V=M\cap T$ holds for some $\pi$-Hall subgroup $T$ of $G$. Since $V=M\cap T\unlhd T$, it follows that $T\leq K$, whence $K\in E_\pi$.

Suppose,  $K<G$. By induction, $K=(K\cap A)N_K(H)$ for some $H\in\Hall_\pi(K\cap A)$. Notice that $|A: (K\cap A)|$ divides $|A:(T\cap A)|$, since $T\leq K$. Now  $(T\cap A)\in\Hall_\pi(A)$, hence $|A: (K\cap A)|$ is a $\pi'$-number and $\Hall_\pi(K\cap A)\subseteq\Hall_\pi(A)$. In particular, $H\in\Hall_\pi(A)$.  We have
$$
G=MK=M(K\cap A)N_K(H)\leq AN_G(H).
$$

Now assume that $K=G$, i.e. $V\unlhd G$. By induction  $$G/M= (A/M)N_{{G/M}}({X/M})=AN_G(X)/M$$ for some ${X/M}\in \Hall_\pi({A/M})$ and thus  $G=AN_G(X)$.

Suppose that $V\ne 1$, i.e. $M=V$. Then $X\in\Hall_\pi(A)$, and the conclusion of the theorem holds.

Consider the remaining case  $V=1$. In other words, consider the case, where $M$ is a $\pi'$-group. By the Schur--Zassenhaus Theorem 
\cite[Theorem~18.1]{Aschbacher}, $X$ possesses a $\pi$-Hall subgroup $H$, moreover $H\in\Hall_\pi(A)$ and $X=HM$. Let $g\in G$. 
In view of the identity $G=AN_G(X)$ there exists  $a\in A$ such that $X^g=X^a$. Now  $H^g$ and $ H^a$ are $\pi$-Hall subgroups of  $X^g$ and, 
by Schur--Zassenhaus theorem,  $H^g$ and $H^a$ are conjugate in $X^g\leq A$. Hence,  $H$ and $H^g$ are conjugate in $A$, the class $H^A$ is $G$-invariant and $G=AN_G(H)$.

 In order to complete the proof, it remains to show that if   $G\in E_\pi$, $A\unlhd G$ and $H\in\Hall_\pi(A)$ is such that  $G=AN_G(H)$, then $N_G(H)\in E_\pi$ and $\Hall_\pi(N_G(H))\subseteq\Hall_\pi(G)$.

Notice that $$
N_G(H)/N_A(H)\cong AN_G(H)/A=G/A\in E_\pi
$$
by Lemma  \ref{base}. Lemma \ref{HallExist} implies that a complete preimage of a  $\pi$-Hall subgroup of $N_G(H)/N_A(H)$ possesses a $\pi$-Hall subgroup  $Y$ such that $Y\cap N_A(H)=H$. Notice that  $|A:(A\cap Y)|$ divides $|A:H|$ and so it is a  $\pi'$-number. Since $|N_G(H):YA|$ is a  $\pi'$-number, $$
|N_G(H):Y|=|N_G(H):YA||YA:Y|=|N_G(H):YA||A:(A\cap Y)|
$$
is a $\pi'$-number as well. Therefore  $Y\in \Hall_\pi(N_G(H))$ and $N_G(H)\in E_\pi$. Now   $$|G:N_G(H)|= |AN_G(H):N_G(H)|=|A:N_A(H)|$$
divides  $|A:H|$ and so is a  $\pi'$-number. Whence $$\Hall_\pi(N_G(H))\subseteq\Hall_\pi(G),$$
and the theorem follows. \qed
\medskip

{\it Proof of  Corollary}~\ref{HN}. Let $G\in C_\pi$, $A\unlhd G$, and $H\in \Hall_\pi(G)$. We need to show that $HA\in C_\pi$.

Let  $$\Gamma=\{K\cap A\mid K\in \Hall_\pi(G)\}.$$ In view of $G\in C_\pi$, $G$ acts transitively on~$\Gamma$ by conjugations.

By Theorem \ref{MainTheorem}, there exists   $X\in \Hall_\pi(A)$ such that $G=AN_G(X)$ (equivalently, $X^G=X^A$), $N_G(X)\in E_\pi$, and $\Hall_\pi(N_X(K))\subseteq \Hall_\pi(G)$. Choose  arbitrary $Y\in\Hall_\pi(N_G(X))$. It is easy to see that   $Y\cap A=X$ and so   $X\in \Gamma$. By the transitivity of  $G$ on $\Gamma$ we have $\Gamma = X^G=X^A$, and this identity implies the transitivity of  $A$ on~$\Gamma$.

Now choose arbitrary  $K\in \Hall_\pi(HA)\subseteq\Hall_\pi(G)$. Since $H\cap A\in \Gamma$ and $K\cap A\in \Gamma$, by the transitivity of $A$ on $\Gamma$,
there exists  $a\in A\leq HA$ such that  $K\cap A=H^a\cap A$. The group $N_{HA}(H^a\cap A)$ is $\pi$-separable, since
all sections of the normal series  $$N_{HA}(H^a\cap A)\unrhd N_{A}(H^a\cap A)\unrhd H^a\cap A\unrhd 1$$ are either $\pi$- or  $\pi'$-groups.
Hence $N_{HA}(H^a\cap A)\in C_\pi$ by Lemma \ref{base1}. Now   $H^a,K\in\Hall_\pi(N_{HA}(H^a\cap A))$, so there exists   $x\in N_{HA}(H^a\cap A)$ such that $K=H^{ax}$.
Thus we obtain that arbitrary  $K\in\Hall_\pi(HA)$ is conjugate with $H$ in  $HA$, whence $HA\in C_\pi$. \qed
\medskip

{\it Proof of  Corollary}~\ref{E_piCrit}. Let  $A\unlhd G$.

Suppose,  $A\in E_\pi$, $G/A\in E_\pi$, and  $H^A=H^G$ for some $H\in \Hall_\pi(A)$. Let $M$ be the full preimage of a $\pi$-Hall subgroup of $G/A$. Then $M\in E_\pi$ in view of
Lemma~\ref{HallExist} and, moreover, $\Hall_\pi(M)\subseteq\Hall_\pi(G)$, since $|G:M|=|G/A:M/A|$ is a $\pi$-number. Thus, $G\in E_\pi$.

The converse statement is a straight consequence of Lemma~\ref{base} and Theorem~\ref{MainTheorem}. \qed
\medskip

{\it Proof of  Corollary}~\ref{ExInvHallSubg}. Let $G\in E_\pi$ and $A\leq \Aut(G)$ with $(|G|,|A|)=1$. Denote by $G^*$ the natural semidirect product of $G$ and $A$. Then $G^*\in E_\pi$
and, by Theorem~\ref{MainTheorem},  there exists $H\in \Hall_\pi(G)$ such that $G^*=GN_{G^*}(H)$. Furthermore, $$N_{G^*}(H)/N_{G}(H)\cong G^*/G\cong A.$$ Hence $(|N_{G^*}(H)/N_{G}(H)|,|N_{G}(H)|)=1$.
In view of the Schur--Zassenhaus theorem \cite[Theorem~18.1]{Aschbacher}, there is  $B\leq N_{G^*}(H)$ isomorphic to $A$ and $A=B^x$ for some $x\in G^*$. Now $$A=B^x\leq  N_{G^*}(H)^x= N_{G^*}(H^x)$$
and $H^x$ is an $A$-invariant $\pi$-Hall subgroup of $G$.\qed

\bigskip


Authors:

REVIN Danila Olegovitch,

Sobolev Institute of Mathematics, pr. Acad. Koptyg, 4,

Novosibirsk, Russia, 630090;

Novosibirsk State University, ul. Pirogova, 2,

Novosibirsk, Russia, 630090.

e-mail: revin@math.nsc.ru

VDOVIN Evgeny Petrovitch,

Sobolev Institute of Mathematics, pr. Acad. Koptyg, 4,

Novosibirsk, Russia, 630090;

Novosibirsk State University, ul. Pirogova, 2,

Novosibirsk, Russia, 630090.

e-mail: vdovin@math.nsc.ru

\end{document}